\documentclass[onefignum,onetabnum]{siamonline171218}

\usepackage{amsfonts}
\usepackage{amsmath}
\usepackage{lmodern}
\usepackage{mathtools}
\usepackage[T1]{fontenc}
\usepackage{graphicx}
\usepackage{epstopdf}
\usepackage{algorithm}
\usepackage{algpseudocode}
\usepackage[T1]{fontenc}
\usepackage{tikz,caption}
\usepackage{subcaption}
\usepackage{pgfplots}
\pgfplotsset{compat=newest}
\pgfplotsset{plot coordinates/math parser=false}
\usepackage{url}
\usepackage{color}
\usepackage{cleveref}

\usetikzlibrary{plotmarks} 
\usetikzlibrary{external}
\usepgfplotslibrary{external} 
\tikzset{external/mode=graphics if exists}
 
\newlength\figureheight
\newlength\figurewidth 

\ifpdf
  \DeclareGraphicsExtensions{.eps,.pdf,.png,.jpg}
\else
  \DeclareGraphicsExtensions{.eps}
\fi

\usepackage{xcolor}

\newcommand{\R}{\mathbb{R}}

\newcommand{\hk}{^{(d)}}
\newcommand{\D}{\,\mathrm{d}}

\def\matlab{{{\sc matlab}}}

\numberwithin{theorem}{section}

\newcommand{\TheTitle}{A low-rank tensor method to reconstruct sparse initial states for PDEs with Isogeometric Analysis} 
\newcommand{\ShortTitle}{Low-rank method for sparsity inducing Priors} 
\newcommand{\TheAuthors}{Alexandra B\"unger, and Martin Stoll}

\headers{\ShortTitle}{\TheAuthors}

\title{{\TheTitle}}

\author{Alexandra B\"unger\thanks{Technische Universit\"at Chemnitz, Department of Mathematics, Chair of Scientific Computing, 09107 Chemnitz, Germany, \email{alexandra.buenger@mathematik.tu-chemnitz.de}}
  \and
  Martin Stoll\thanks{Technische Universit\"at Chemnitz, Department of Mathematics, Chair of Scientific Computing, 09107 Chemnitz, Germany, \email{martin.stoll@mathematik.tu-chemnitz.de}}
}
\usepackage{amsopn}
\DeclareMathOperator{\diag}{diag}

\ifpdf
\hypersetup{
  pdftitle={\TheTitle},
  pdfauthor={\TheAuthors}
}
\fi

\def\phys{} 

\definecolor{grey}{rgb}{0.5,0.5,0.5}
\definecolor{darkgreen}{rgb}{0,0.55,0}

\begin{document}

\maketitle

\begin{abstract}
When working with PDEs the reconstruction of a previous state often proves difficult. Good prior knowledge and fast computational methods are crucial to build a working reconstruction. We want to identify the heat sources on a three dimensional domain from later measurements under the assumption of small, distinct sources, such as hot chippings from a milling tool. This leads us to the need for a Prior reflecting this a priori information. Sparsity-inducing hyperpriors have proven useful for similar problems with sparse signal or image reconstruction. We combine the method of using a hierarchical Bayesian model with gamma hyperpriors to promote sparsity with low-rank computations for PDE systems in tensor train format.
\end{abstract}



\section{Motivation}
The reconstruction of sparse sources from noisy observation of a PDE model plays a significant role in many applications, varying from identifying heat sources \cite{knapik2013} to finding the origin of a tumor growth \cite{subramanian2020}.
Incorporating prior information about the sparsity of an unknown source into its recovery from noisy data arose as long as 40 years ago with first successful applications in geophysics, astrophysics, and ultrasonic imaging \cite{Donoho1992}. 

The general idea of using $\ell_p$ Tikhonov regularization techniques has since become highly popular to promote sparsity and has been extensively researched (cf. \cite{Grasmair2008,Jin2017} and references therein), mainly focusing on $1 \leq p \leq 2$ due to their global convexity properties. It has been shown that an $\ell_1$ penalty can provide accurate sparse reconstructions \cite{Candes2006} or even the sparsest possible solution under some additional assumptions \cite{Donoho2006}. Stronger sparsity promotion with $p < 1$ comes at the expense of global convexity and is therefore still an active research topic \cite{Zarzer2009}. 

In the Bayesian framework, sparsity promotion can be efficiently modeled by considering a component-wise Gaussian prior model where the variances themselves are modeled as random variables underlying a sparsity promoting \textit{hyperprior} model \cite{Calvetti2019}.
In this paper we will examine the application of a hybrid method using such hypermodels with gamma hyperpriors presented in \cite{calvetti2020}. The algorithm strongly promotes sparsity while retaining global convergence by switching between two models where the first is globally convex to drive the solution close to a unique minimizer and the other is more greedy towards stronger sparsity while suffering from local minima. 

Let us assume a PDE model with Dirichlet boundary conditions on a 2D or 3D domain $\Omega$ as
\begin{align}
\frac{\partial}{\partial t} y(x,t) &= \mathcal{L} y(x,t) &&\mbox{ on } \Omega \times (0,T], \label{eq:PDE1}\\
y(x,t=0) &= u_0(x) &&\mbox{ on } \Omega \times 0, \label{eq:PDE2} \\
y(x,t) &=y_0(x) &&\mbox{ on } \partial \Omega \times (0,T], \label{eq:PDE3}
\end{align}
with a PDE operator $\mathcal{L}$ such as $\mathcal{L}y(x,t) = \Delta y(x,t)$ for heat conduction. Here, $y(x,t)$ is the state at time $t$ and $u_0(x,t)$ the respective initial state at time $t=0$. Both are time-dependent functions on the domain $\Omega$ and time-frame $(0,T]$ with end time $T$.

We consider being given a noisy observation at the end time $T$ stemming from the evolution via the PDE model \cref{eq:PDE1,eq:PDE2,eq:PDE3},
\begin{equation}
z(x) = y(x,T) + e.
\end{equation}
Our goal is to recover the initial state under the presence of additive noise $e$, e.g., measurement errors, and the prior assumption of $u_0$ being sparse.

In this paper, we will study a Bayesian approach to the problem of reconstructing $u_0$ under the assumption of Gaussian additive noise and a sparsity inducing prior formulation. We  make use of the methods developed in \cite{calvetti2020} and tailor them to our model in \cref{chapter:prior}. The original algorithm requires a large number of matrix-vector multiplications with the parameter-to-observable map, which in our case includes a complete solution of the governing PDE and thus can be infeasibly large and dense.

The PDE model can be discretized in space and time into a linear equation system using some time-stepping method and a Galerkin type discretization, which we will introduce in \cref{chapter:discretization}. We introduce the space discretization using an Isogeometric Analysis (IGA) scheme in \cref{subsection:iga} and the full space-time discretization in \cref{subsection:discretePTO}. As the resulting discretization would be large and dense, we propose a low-rank approximation method for the IGA discretization in \cref{chapter:LowRank} and a suitable solver tailored to working with low-rank tensor structured data in \cref{chapter:amen} to perform matrix-vector multiplications or solving the arising KKT system without computing the matrices directly. The resulting numerical scheme profits from small storage requirements and robustness allowing to compute even very large-scale problems, which would prove infeasible with traditional methods. In \cref{chapter:examples} we provide numerical examples with an exemplary 2D problem and 3D models to show the effectiveness and robustness of the resulting method comparing different hyperprior setups. 

\section{Sparsity inducing hyper prior} \label{chapter:prior}

We will first briefly review the method for sparsity inducing hyperpriors as the local hybrid IAS (iterative alternating sequential) algorithm introduced in \cite{calvetti2020} and later extended in \cite{calvetti2020:2} from a Bayesian point of view.

We consider the linear observation model with additive Gaussian noise resulting from the PDE problem with uncertain initial values $u_0$,
\begin{equation}
	z = Au_0 + e, \quad e \sim \mathcal{N}(0,\Sigma),
\end{equation}
where $A \in \R^{n \times m}$ is the parameter-to-observable map, which transfers the uncertain inputs $u_0$ to the observable outputs $y$, and $\Sigma \in \R^{m \times m}$ is a symmetric positive semidefinite covariance matrix. Without loss of generality, we assume the noise to be white, i.e. $\Sigma = I$. This gives the likelihood probability density function (pdf)
\begin{equation}
	\pi_{z|u_0}(z|u_0) \propto \exp \left ( - \frac{1}{2} \| A u_0 - z\|^2 \right ).
\end{equation}
We believe the unknown $u_0$ to be sparse. To incorporate this into the prior model, we consider a component-wise Gaussian prior model,
\begin{equation}
u_0 \sim \mathcal{N}(0,D_\theta), \, D_\theta = \diag(\theta_1, \hdots, \theta_n) \in \R^{n\times n},
\end{equation}
where the variances $\theta$ are not known. Thus, $\theta$ itself are random variables where smaller variances promote values closer to zero. The a priori belief about $\theta$ is thus incorporated into a \textit{hyperprior} pdf $\pi_\theta(\theta)$ leading to a hierarchical prior model. The conditional prior density of $u_0$ with given $\theta$ now has the form
\begin{equation}
\pi_{u_0 | \theta} (u_0 |\theta) \propto \frac{1}{\prod_{j=1}^n \sqrt{\theta_j}} \exp \left (-\frac{1}{2}\|D_\theta^{-1/2}u_0\|^2 \right) = \exp \left ( - \frac{1}{2}\|D_\theta^{-1/2}u_0 \|^2 - \frac{1}{2} \sum_{j=1}^n \log \theta_j \right).
\end{equation}
Now we need to estimate $\theta$ as well as $u_0$ and chose an appropriate prior density for its representation. One class of hyperpriors $\pi_\theta$ to promote the sparsity of our solution is that of generalized gamma distributions,
\begin{equation} \label{eq:gamma}
\pi_\theta(\theta) = \pi_\theta(\theta|r, \beta, \vartheta) = \frac{|r|^n}{\Gamma(\beta)^n} \prod_{j=1}^n \frac{1}{\vartheta_j} \left (\frac{\theta_j}{\vartheta_j} \right )^{r\beta-1} \exp \left ( - \left ( \frac{\theta_j}{\vartheta_j} \right )^r \right),
\end{equation}
with $r \in R \backslash \{0\}, \beta > 0, \vartheta > 0$. This family of priors allows large outliers in $\theta$ while overall favoring small values \cite{calvetti2020}.

With this we not only estimate $x_0$ but $\theta$ as well, based on their joint posterior distribution,
\begin{equation}
\pi_{u_0,\theta} (u_0,\theta | y) \propto \pi_{u_0|\theta}(u_0|\theta) \pi_{\theta}(\theta) \pi_{z|u_0}(z|u_0). \label{eq:posterior_pdf}
\end{equation}
To find a suitable sparse reconstruction for $u_0$, we need to compute the maximum a posteriori (MAP) estimate, which is the minimizer of the negative logarithm of the posterior pdf \eqref{eq:posterior_pdf},
\begin{align}
(u_0^*, \theta^*) &= \arg \min_{u_0, \theta} \mathcal{F}(u_0,\theta | r, \beta, \vartheta) \\
&= \arg \min_{u_0, \theta} \frac{1}{2}\|z - Au_0\|^2 + \frac{1}{2} \sum_{j=1}^n \frac{u_{0,j}^2}{\theta_j} - \left ( r \beta - \frac{3}{2} \right ) \sum_{j=1}^n \log \frac{\theta_j}{\vartheta_j} + \sum_{j=1}^n \left (\frac{\theta_j}{\vartheta_j} \right )^r. \label{eq:minimization}
\end{align}
Note that the objective function $\mathcal{F}(u_0, \theta)$ consists of four terms of which one depends on $u_0$ only, two depend only on $\theta$ and one is dependent on both variables. This gives rise to an hybrid algorithm alternating between updates of $u_0$ and $\theta$.

The authors of \cite{calvetti2020} suggest using an Iterative Alternating Sequential algorithm (IAS) for this computation. Here, each iteration consists of two updates,
\begin{equation}
u^i, \theta^i \rightarrow u^{i+1} \rightarrow \theta_{i+1},
\end{equation}
where we first fix $\theta$ to update $u^{i+1}$ and afterwards fix $u^{i+1}$ to update $\theta$,
\begin{equation}
u^{i+1} = \arg \min_u \{\mathcal{F}(u,\theta^i)\}, \, \theta^{i+1} = \arg \min_\theta \{\mathcal{F}(u^{i+1}, \theta)\}.
\end{equation}
Both minimizations have an exact condition for the minimizer and are relatively simple to compute. The update of $u$ results in the solution of a quadratic minimization problem, i.e.,
\begin{equation}
u^{i+1} = \arg \min_u \|Au - z\|^2 + \| D_\theta^{-1/2}u\|^2. \label{eq:LS}
\end{equation}
To approximate the solution, we can solve \cref{eq:LS} as it is or solve a reduced least squares problem. For this we introduce a change of variables as $w = D^{-1/2}_\theta u$ and find the least squares solution of
\begin{equation}
A D^{1/2}_\theta w = z. \label{eq:reduced}
\end{equation}
This can be easily solved, e.g. by a CGLS (Conjugate Gradient Least Squares) algorithm, as suggested in \cite{calvetti2020:2}. We will later discuss the structure of $A$ for our PDE model and propose a low rank method to efficiently solve this minimization problem.

To then update $\theta$ we can compute its components $\theta_j$ independently from each other. For each component the first order optimality conditions applied to \eqref{eq:minimization} read
\begin{equation}
0 = \frac{\partial \mathcal{F}}{\partial \theta_j} = - \frac{1}{2} \frac{u_j^2}{\theta_j^2} - \left ( r \beta - \frac{3}{2} \right ) \frac{1}{\theta_j} + r \frac{\theta_j^{r-1}}{\vartheta_j^r}.
\end{equation}
To numerically solve this equation, \cite{calvetti2020} suggests the change of variables $\theta_j = \vartheta_j \xi_j$, $u_j = \sqrt{\vartheta_j}z_j$. Then, we can write $\xi_j = \varphi(|z_j|)$ and we get the initial value problem
\begin{equation}
\varphi'(z) = \frac{2 z \varphi(z)}{2r^2\varphi(z)^{r+1} + z^2}, \, \varphi(0) = \left ( \frac{\eta}{r}\right )^{1/r}.
\end{equation}
From this the update for $\theta_j$ can be computed with some numerical time integrator.

Alternating between updating $u$ and updating $\theta$ gives rise to the standard IAS algorithm. 

\subsection{Hyperprior parameters and hierarchical IAS algorithm} \label{chapter:hyperpriors}

The effectiveness of the alternating scheme depends highly on the choice of parameters for the hyperprior, $r$, $\beta$ and $\vartheta$. A detailed analysis of these parameters is available in \cite{calvetti2020}, which we want to summarize now. We additionally analyze the choice of parameters in our numerical experiments in \Cref{chapter:examples}.

The choice of $r$ and $\beta$ affects the sparsity of the solution and determines the convexity of the objective function while $\vartheta$ can be seen as a sensitivity scaling where $\vartheta_j = \frac{C}{\|Ae_j\|^2}$ with some constant $C > 0$ (cf. \cite{calvetti2020:2}). The constant $C$ should depend on the signal-to-noise ratio (SNR).

Calvetti et al. propose two modifications to the general IAS algorithm to update the hypermodel component-wise and ensure convexity. When $r > 1$ and $\beta$ follows $r\beta>3/2$ the objective function is globally convex. Additionally, the objective function is convex for $0<r<1$, $r\beta >3/2$ or $r<0$ and $\beta >0$ when $\theta_j < \vartheta_j(\frac{\eta}{r|r-1|})$. This can be exploited by choosing two hypermodels $M_1$ and $M_2$, where $M_1$ satisfies the global convexity conditions and $M_2$ is the desired sparsity promoting model. 

\textbf{Global hybrid IAS.} The first method proposed starts the iteration with model $M_1$ and caries out a switch to model $M_2$ for all parameters at once after a certain number of iterations. The first iterations drive the objective function towards the global minimum and afterwards we trade global convexity for stronger sparsity promotion. The algorithm is outlined in \cref{alg:global}.

\textbf{Local hybrid IAS.} The second proposition starts with model $M_2$ as well but updates the models locally by switching individual components of $\theta_j$ to the other model $M_2$ when $\theta_j$ is below the threshold for convexity. To ensure preservation of convexity an additional bound constraint is applied to components, which were switched to model $M_2$. \Cref{alg:local} illustrates the iteration for this method.

We will consider both methods in our experiments for a gamma ($r=0.5$) and inverse gamma ($r=-1$) hyperprior model. 

 \begin{algorithm}
 	\caption{Global hybrid IAS} \label{alg:global}
 	\begin{algorithmic}[1]
 		\State Input: noisy $z$, parameter-to-observable map $A$, noise covariance matrix $\Sigma$, hypermodel $M_1$ with $(r_1, \beta_1, \vartheta_1)$, hypermodel $M_2$ with $(r_2, \beta_2, \vartheta_2)$, switch point $i_s > 0$, 
 		\State Set $\theta_0 = \vartheta_1$
 		\For{$i=1,\hdots,\text{maxIter}$ or until convergence}
	 		\State Solve $\min_{u_i}  \|A u_i -z \|^2 + \|D_\theta^{-1/2} u_i \|^2$
	 		\If{$i<i_s$}
		 		\State Update $\theta$ with $u_i$ and parameters from $M_1$
	 		\Else
		 		\State Update $\theta$ with $u_i$ and parameters from $M_2$
	 		\EndIf
	 	\EndFor
		\State Output: estimated initial values $u_0 = u_i$ and variance $\theta$
 	\end{algorithmic}
 \end{algorithm}
 
  \begin{algorithm}
  	\caption{Local hybrid IAS} \label{alg:local}
  	\begin{algorithmic}[1]
  		\State Input: noisy $z$, parameter-to-observable map $A$, noise covariance matrix $\Sigma$, hypermodel $M_1$ with $(r_1, \beta_1, \vartheta_1)$, hypermodel $M_2$ with $(r_2, \beta_2, \vartheta_2)$
  		\State Set $\theta_0 = \vartheta_1$, $I = \emptyset$
  		\For{$i=1,\hdots,\text{maxIter}$ or until convergence}
	  		\State Solve $\min_{u_i}  \|A u_i -z \|^2 + \|D_\theta^{-1/2} u_i \|^2$
	  		\State Apply bound constraint to $u_{i,j}$ with $j \in I$
	  		\For{$j=1,\hdots,n$}
		  		\If{$\theta_j \geq \vartheta_j(\frac{\eta}{r|r-1|})$}
			  		\State update $\theta_j$ with model $M_1$
			  	\Else
				  	\State update $\theta_j$ with model $M_2$
				  	\State $i = I \cup \{j\}$
				 \EndIf
			\EndFor
		\EndFor
  		\State Output: estimated initial values $u_0 = u_i$ and variance $\theta$
  	\end{algorithmic}
  \end{algorithm}

\section{Space-time discretization of $A$} \label{chapter:discretization}

To use the method described in \cref{chapter:prior} we need a representation of the parameter-to-observable map $A$. For our PDE model \cref{eq:PDE1,eq:PDE2,eq:PDE3}, we need an appropriate discretization to build this mapping. We will now introduce Isogeometric Analysis for the space discretization in \cref{subsection:iga} and the time-discretization with an implicit Euler scheme in \cref{subsection:discretePTO}. 

For the space discretization a Finite Element Method (FEM) \cite{FEM} is the most prominently used choice and has seen a wide range of applications in both engineering and science \cite{singiresu2018}. But recently a new approach called Isogeometric Analysis (IGA) \cite{hughes2009} has gained popularity. This approach arose from the desire to integrate numerical analysis directly into CAD (computer aided design) tools \cite{hughes2009}. It abolishes the need for time consuming meshing steps and preprocessing needed in classical Finite Element Analysis \cite{zhu2018} in addition to working on exact geometric representations instead of an approximation.

Furthermore, this method admits a tensor-product based structure which gives rise to a very effective low-rank method developed in \cite{angelos1} and refined in \cite{angelos2, angelos3}. This has been shown to be very efficient in combination with tensor train calculations \cite{osel-tt-2011} in our recent work on PDE-constrained optimal control problems \cite{buenger2020}.

\subsection{Isogeometric Analysis discretization} \label{subsection:iga}

Isogeometric Analysis benefits from the use of high order spline functions, such as B-splines and NURBS, as their basis. A geometric domain is represented exactly using a set of such splines, which are subsequently used to build a solution space for solving a PDE problem \cite{CAD}. In this paper we will focus on the use of B-splines but note that a generalization to other popular spline spaces, such as NURBS (Non-uniform rational B-splines) \cite{NURBS,NURBS2}, is possible. The notation and derivation in this section as well as \cref{chapter:LowRank} follow our work in \cite{buenger2020}.

A set of $n$ B-splines is uniquely defined by choosing a vector $\xi= \{\hat{x}_1, \hdots \hat{x}_{n+p+1}\}$, called the open knot vector, with
\begin{equation} 
0= \hat{x}_1 = \hdots = \hat{x}_{p+1} < \hat{x}_{p+2} \leq \hdots \leq \hat{x}_n < \hat{x}_{n+1} = \hdots = \hat{x}_{n+p+1} = 1, 
\label{equation:knotVector}
\end{equation}
 and a degree $p$, where the first and last knot are repeated $p+1$ times and for all other knots duplicate knots are allowed up to multiplicity $p$. The parameter $n$ determines the number of resulting B-splines $\beta_{i,p}$ with $i=1,\hdots,n$. 

\begin{figure}[H]
\centering
\includegraphics[width = 0.8\textwidth]{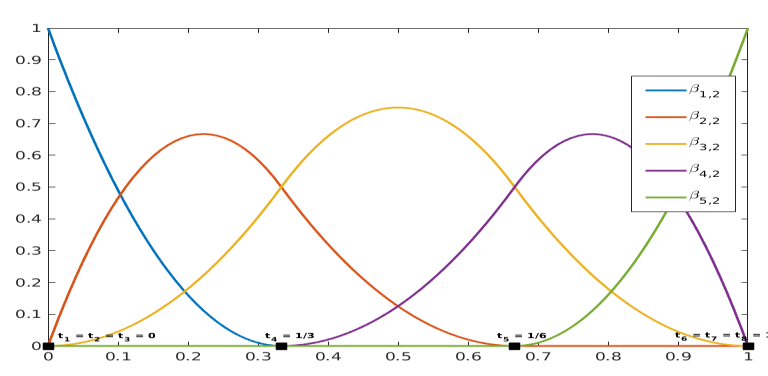}
\caption{B-spline space with degree $p=2$ and knot vector $\xi = [0,0,1/3,1/2,2/3,1,1]$.} \label{fig:spline}
\end{figure}

For each knot vector $\xi$ as in \cref{equation:knotVector}, the according  B-splines $\beta_{i,p}$ of degree $p$ with $i = 1,\hdots, n$ are uniquely defined by the recursion
\begin{align}
\beta_{i,0}(\hat{x}) &= \begin{cases} 1 & \mbox{if } \hat{x}_i \leq \hat{x} < \hat{x}_{i+1}, \\ 0 & \mbox{otherwise}, \end{cases} \\
\beta_{i,j}(\hat{x}) &= \frac{\hat{x}-\hat{x}_i}{\hat{x}_{i+j} - \hat{x}_{i}} \beta_{i, j-1}(\hat{x}) + \frac{ \hat{x}_{i+j+1} - \hat{x}}{\hat{x}_{i+j+1} - \hat{x}_{i+1}} \beta_{i+1, j-1}(\hat{x}),
\end{align}
where $j = 1,2,\hdots, p$ and $i = 1, \hdots, n$. Each resulting B-spline $\beta_{i,p}$ has the local support $[\hat{x}_i,\hat{x}_{i+p+1}]$, see \cref{fig:spline} for an example. We use $\mathbb{S}_\xi^p$ to denote the spline space spanned by the B-splines with degree $p$ and knot vector $\xi$,
\begin{equation}
\mathbb{S}_\xi^p = \text{span} \{ \beta_{1,p}, \hdots, \beta_{n,p} \}.
\end{equation} 

To construct a space of $D$-dimensional basis functions $\mathbb{S}^D$, for $d=1,\hdots,D$ we define one-dimensional spline spaces $\mathbb{S}_{\xi_d}^{p_d}$ and vector valued functions $B\hk(\hat{x}\hk)$, as the vectors holding all basis functions of dimension $d$,
\begin{equation}
B\hk (\hat{x}\hk )  = [ \beta\hk_{1}(\hat{x}\hk ), \hdots, \beta\hk_{n_d}(\hat{x}\hk ) ], \label{eq:Bd}
\end{equation}
with one-dimensional variables $\hat{x}\hk$.

We combine these spline spaces via tensor product to get a D-dimensional spline space $\mathbb{S}^D =  \mathbb{S}_{\xi_1}^ {p_1} \otimes \cdots \otimes \mathbb{S}_{\xi_D}^{p_D}$ with an order $D$ tensor of size $(n_1,\hdots, n_D)$ collecting all basis functions as
\begin{equation}
B(\hat{x}) = \bigotimes_{d=1}^D B\hk(\hat{x}\hk). \label{eq:Basis}
\end{equation}

The tensor  in \cref{eq:Basis} then is a function $B(\hat{x}) : \R^D \rightarrow \R^{n_1 \times \hdots \times n_D}$ with variables $\hat{x} = (\hat{x}^{(1)}, \hdots, \hat{x}^{(D)})^T$ and its elements are the $D$-dimensional basis functions 
\begin{equation}
\beta_\mathbf{i} (\hat{x}) = \prod_{d=1}^D \beta_{i_d}\hk (\hat{x}\hk),
\end{equation}
with multi-index $\mathbf{i} \in I =  \{ (i_1, \hdots, i_D) \, | \, i_d = 1,\hdots, n_d ,\, d=1,\hdots, D\}$.

Given such a basis $\mathbb{S}^D$, we define a B-spline geometry mapping $G:\hat{\Omega} \rightarrow \Omega$ from the $D$-dimensional unit cube $\hat{\Omega}:=[0,1]^D$ onto an arbitrary geometric shape $\Omega \subset \R^D$ as
\begin{equation} \label{equation:geometryMapping}
G(\hat{x}) = \sum_{\mathbf{i} \in I} C_\mathbf{i} \beta_{\mathbf{i}}(\hat{x}) = C:B(\hat{x}),
\end{equation}
with control points $C_{\mathbf{i}} \in \R^D$. Here $C:B(\hat{x})$ denotes the Frobenius product of the two tensors, $C \in \R^{D \times n_1 \times \hdots \times n_D}$ holding all the control points and $B$ from \cref{eq:Basis}.

Now that we have a spline representation of the geometry $\Omega$, we use the same splines to construct discrete functions $u_h \in V_h \subset H^1(\Omega)$ as approximations to the solutions $u \in H^1(\Omega)$ of the PDE problem in \cref{eq:PDE1,eq:PDE2,eq:PDE3}.

The isogeometric solution space $V_h$ is built as
\begin{equation} 
V_h = \mbox{span}\{\beta_\mathbf{i} \circ G^{-1} \, \, : \, \, \mathbf{i} \in I \}, \label{eq:Vh}
\end{equation}
with an index set $I$ such that $\beta_{\mathbf{i}}$ are the elements of $\mathbb{S}^D$.
The functions in $V_h$ are linear combinations of the basis functions with coefficients $u_\mathbf{i}$,
\begin{equation} 
u_h = \sum_{\mathbf{i} \in I} u_{\mathbf{i}} ( \beta_{\mathbf{i}} \circ G^{-1}).
\label{equation:coefficientSet}
\end{equation}

We use this space for the spatial Galerkin discretization of the PDE system \cref{eq:PDE1,eq:PDE2,eq:PDE3}. Exemplary, for the boundary value problem $u(x) = f(x)$ in $\Omega$, the weak formulation results in the bilinear form
\begin{equation}
a_m(u,v) = \langle u, v \rangle_2 = \int_\Omega uv \D x,
\end{equation}
called the mass term. Discretization with \cref{eq:Vh} results in the discrete mass term\begin{equation}
a_{m,h}(u_h, v_h) = \int_\Omega u_h(x) v_h(x) \D x = \int_{\hat{\Omega}} \sum_{\mathbf{i} \in I} u_\mathbf{i} \beta_\mathbf{i}(\hat{x}) \sum_{\mathbf{j} \in I} v_\mathbf{j} \beta_\mathbf{j}(\hat{x}) \omega(\hat{x}) \D \hat{x},
\end{equation}
with $\omega(\hat{x}) = | \det \nabla G(\hat{x})| $. The basis functions $\beta_{\mathbf{j}}$ are linearly independent. Thus, similarly to FEM, we can rewrite the bilinear forms as a matrix-vector product $Au$ with vectorization of the coefficient set $u_\mathbf{i}$, $\mathbf{i} \in \mathcal{I}$, where $A$ is realized as a mass matrix $M$ with elements
\begin{equation} \label{equation:massMatrix}
M_{\mathbf{i},\mathbf{j}} = \int_{\hat{\Omega}} \beta_\mathbf{i} \beta_\mathbf{j} \omega \D \hat{x}.
\end{equation}

The same strategy is applied to any PDE operator $\mathcal{L}u(x)$, e.g. $\mathcal{L}u(x) = - \Delta u(x)$ results in the stiffness term
\begin{align}
a_s(u,v) &=  -\int_\Omega  (\Delta u) v \D x = \int_\Omega \nabla u \cdot \nabla v \D x, \label{equation:stiffness}, \\
a_{s,h}(u_h,v_h) &= \int_{\Omega} (\nabla u_h(x) ) \cdot \nabla v_h(x) \D x = \int_{\hat{\Omega}}(Q(\hat{x}) \sum_{\mathbf{i} \in I} u_\mathbf{i} \nabla \beta_\mathbf{i}(\hat{x})) \cdot\sum_{\mathbf{j} \in I} v_\mathbf{j} \nabla \beta_\mathbf{j}(\hat{x}) \D \hat{x},
\end{align}
with $Q(\hat{x}) = \big (\nabla G(\hat{x})^{T} \nabla G(\hat{x})\big )^{-1} | \det \nabla G(\hat{x})|$ (cf. \cite{angelos1}) resulting in the stiffness matrix $K$ with elements
\begin{equation} \label{equation:stiffnessMatrix}
K_{\mathbf{i},\mathbf{j}} = \int_{\hat{\Omega}}( Q \nabla \beta_\mathbf{i} )\cdot \nabla \beta_\mathbf{j}  \D \hat{x} = \sum_{k,l=1}^D \int_{\hat{\Omega}} q_{k,l} \frac{\partial}{\partial \hat{x}_l} \beta_\mathbf{i}  \frac{\partial}{\partial \hat{x}_k} \beta_\mathbf{j} \D \hat{x}.
\end{equation}

With these matrices we can now discretize the parameter-to-observable map $A$ in a general way before introducing additional low-rank approximations for $M$ and $K$ exploiting the underlying tensor product structure of the basis functions to reduce the computational effort and storage requirements of our problem. This low-rank approach has been successfully applied before to an optimal control problem in \cite{buenger2020}.

\subsection{Discrete parameter-to-observable map} \label{subsection:discretePTO}

Considering the PDE problem given in \eqref{eq:PDE1} -- \eqref{eq:PDE3}, we now have a space discretization, which we need to combine with a suitable method for time discretization. For this time-stepping method we choose an implicit Euler scheme \cite{euler} and $N_t$ equidistant time steps of length $\tau$ to discretize the time frame $t \in (0,T]$ into steps $t_i = i \tau$ for $i = 1,\hdots, N_t$. With this we get the fully discrete system
\begin{align}
\frac{\phys{M} y_i - \phys{M} y_{i-1}}{\tau} &= - \phys{K} y_i,  &&\mbox{ for } i=2,\hdots, N_t \\
M y_1 &= -\tau K y_1 + M u_0 \\
z &= C y_{N_t} + e.
\end{align}
Here, $y_i$ is the vector of B-spline coefficients for the discretization of the function $y(t_i)$ and $C$ denotes the observation operator (e.g. an identity for a full observation at $N_t$). We will denote the vectors collecting all time steps as e.g. $\mathbf{y} = [y_1, \hdots, y_{N_t}]^T$.

We can write this system of equations as one large equation system,
\begin{align}
\underbrace{\begin{bmatrix} \tau \phys{K} + \phys{M} & & \\- \phys{M}  & \tau \phys{K} + \phys{M} & \\ &  \hspace{-3em}\ddots & \hspace{-3em} \ddots \\  & -\phys{M} & \tau \phys{K} + \phys{M} \end{bmatrix}}_{\mathbf{K}} \begin{bmatrix} y_1 \\ \vdots \\ y_{N_t} \end{bmatrix} &=
\underbrace{\begin{bmatrix} \phys{M} \\ 0 \\ \vdots \\ 0 \end{bmatrix}}_{\mathbf{M}_0} u_0, \\
 z = \underbrace{\begin{bmatrix} 0 & \cdots & 0  & C \end{bmatrix}}_{\mathbf{C}} \begin{bmatrix} y_1 \\ \vdots \\ y_{N_t} \end{bmatrix} &+ e.
\end{align} 
Using this formulation we can denote the equation from $u_0$ to $z$, the parameter-to-observable map $A$, as
\begin{equation} 
z = Au_0 := \mathbf{C} \mathbf{K}^{-1} \mathbf{M}_0 u_0. \label{eq:parameter_to_observables}
\end{equation}
This matrix $A$ typically is dense and large, and even storing the full representation of $\mathbf{K}$ may pose difficulties. Thus, computing $A$ explicitly is usually infeasible. We will present a scheme to efficiently compute matrix-vector products with $A$ without forming any full matrices. For this, let us observe the structure of the arising matrices. We can also expressed them with the following Kronecker products,
\begin{align}
\mathbf{C} &= \begin{bmatrix} 0 & \cdots &0 & 1\end{bmatrix} \otimes C, \label{eq:kronecker1} \\
\mathbf{K} &= \tau I_{N_t} \otimes \phys{K} + II \otimes \phys{M},  \label{eq:kronecker3}\\
\mathbf{M}_0 &= \begin{bmatrix} 1 & 0 & \cdots & 0\end{bmatrix}^T  \otimes \phys{M}, \label{eq:kronecker4}
\end{align}
where $I_{N_t}$ is the $N_t \times N_t$ identity, $II$ is a bidiagonal matrix with 1 on the diagonal and -1 on its  first subdiagonal\footnote{Syntax for \matlab \,: $\mathtt{II = spdiags([-ones(N_t,1),ones(N_t,1)],-1:0,N_t,N_t);}$}.

Note that here $II$ results from the implicit Euler scheme. This scheme can be replaced by a different time stepping method like a Crank-Nicolson scheme \cite{juncosa1957} or others in which case $II$ and the identity $I_{N_t}$ will become different block matrices. Using a different time-stepping method would thus still maintain the general Kronecker product structure of the equations and be applicable for the following steps.

Now that we know the structure of our parameter-to-observable map we can take a closer look at solving the minimization problem in \cref{eq:LS}. In \cite{calvetti2020}, the authors suggest using a CGLS method \cite{CGLS}. This method is most suited if the underlying matrix $A$ is sparse and can be precomputed (cf. \cite{Shewchuk1994}) as it requires a large number of matrix-vector products with $A$ and $A^T$. 

Unfortunately, in the case of the PDE solution \cref{eq:parameter_to_observables} the resulting matrix will most likely be dense and potentially very large. Thus, constructing $A$ in full may not be feasible. In the following chapter we propose a low rank method that can be used to either calculate the matrix-vector products needed for CGLS in a compact and efficient way without assembly of $A$, or to solve the optimization problem \cref{eq:LS} by rearranging it into a constrained minimization problem.

Substituting $A u_0$ in \cref{eq:LS} by $\mathbf{C}y$ and using the PDE as a constraint leads to
\begin{align}
\min_{u_0} \,\|\mathbf{C}y - z\|^2 &+ \|D_{\theta}^{-1/2}u_0\|^2 \label{eq:constrained1} \\
\text{s.t. }\mathbf{K} y &= \mathbf{M}_0 u_0. \label{eq:constrained2}
\end{align}
Formulating its first order optimality conditions gives us the large scale saddle point system,
\begin{align}
\begin{bmatrix} \mathbf{C}^T M \mathbf{C} & 0 & -\mathbf{K}^T \\ 0 & D^{-1/2}_{\theta} M D^{-1/2}_{\theta} & \mathbf{M}_0^T \\ -\mathbf{K} & \mathbf{M}_0 & 0 \end{bmatrix} \begin{bmatrix} y \\ u_0 \\ \lambda \end{bmatrix} = \begin{bmatrix} \mathbf{M}_0 \\ 0 \\0 \end{bmatrix}. \label{eq:saddlepoint}
\end{align}
This type of system has been extensively researched and can generally be solved with e.g. preconditioned iterative methods (cf. \cite{Pearson2012}) or variational formulations \cite{Hinze2005}. But due to its potential size we will focus on a low rank method exploiting the underlying structure of the system. Here, each block has a Kronecker product structure as described in \cref{eq:kronecker1,eq:kronecker1,eq:kronecker1}, which can be exploited efficiently by a low rank in time method presented in \cite{stoll1}. To further reduce the computational complexity we will now take a look back at the structure of the matrices $M$ and $K$ and outline how we can solve the resulting Kronecker product based system with a low rank tensor train method.

\section{Low Rank IGA} \label{chapter:LowRank}
We see that the time discretization in \cref{eq:kronecker1,eq:kronecker3,eq:kronecker4} leads to a tensor product structure in the parameter-to-observable map. Using the method presented in \cite{buenger2020} we can additionally approximate $M$ and $K$ as similar tensor products and exploit the resulting structure with low rank tensor train calculations. This exploit can either be utilized for the matrix-vector products in the CGLS scheme or to solve the saddle point system \cref{eq:saddlepoint}.

During the derivation of the mass and stiffness matrices, we did not pay attention to the tensor product structure of $\mathbb{S}^D$. We can either arrange $M$ and $K$ as matrices or as tensors of size $(\mathbf{n},\mathbf{n}) = (n_1,\hdots,n_D,n_1,\hdots,n_D)$. With this tensor notation the mass and stiffness matrices in a multi-dimensional setting are represented in a compact way. We can write the mass term as a tensor $M$
\begin{equation} \label{equation:massTensor1}
M = \int_{\hat{\Omega}} \omega B \otimes B \D \hat{x} \, \, \, \in \R^{\mathbf n \times \mathbf n},
\end{equation} 
with elements coming from \cref{equation:massMatrix}. The stiffness term can be treated similarly. With the tensor gradient we can write it as a tensor $K$
\begin{align} \label{equation:stiffnessTensor1}
K &= \int_{\hat{\Omega}} [Q \cdot( \nabla \otimes B)] \cdot (\nabla \otimes B) \D \hat{x} \, \, \, \in \R^{\mathbf n \times \mathbf n}.\\
&= \sum_{k,l=1}^D \int_{\hat{\Omega}} q_{k,l} \frac{\partial}{\partial \hat{x}_l} B \otimes \frac{\partial}{\partial \hat{x}_k} B \D \hat{x}
\end{align}
whose elements are of the form in \cref{equation:stiffnessMatrix}.
The associated mass and stiffness matrices are obtained by reordering the indices since the elements of the mass and stiffness tensors match the elements of the matrices. 

Except for the D-variate weight functions $\omega$ and $Q$, the entries of the mass and stiffness tensors \cref{equation:massTensor1,equation:stiffnessTensor1} are the product of univariate B-splines. The scalar $\omega(\hat{x}) = |\mbox{det } \nabla G(\hat{x})|$ and the matrix $Q(\hat{x}) = (\nabla G(\hat{x})^{-1}) (\nabla G(\hat{x}))^{-T} \omega(t) \in \R^{D\times D}$ are determined by the geometry mapping. As Mantzaflaris et al. suggest in \cite{angelos1}, we can approximate these weight functions via interpolation by some combination of univariate functions,
\begin{equation}
\omega(\hat{x}) \approx \omega_1(\hat{x}^{(1)}) \cdots \omega_D(\hat{x}^{(D)}).
\end{equation}
The integrands then are separable into products of univariate integrals. To further reduce the computation time and storage requirements of the mass and stiffness matrix calculation, the resulting interpolating function is approximated with low rank methods giving low rank approximations of the system matrices \cite{angelos2,angelos1}.

To do so, we interpolate the weight functions by a combination of univariate B-splines of higher order, denoted by the spline space $\tilde{\mathbb{S}}^D$ with suitable knot vectors $\tilde{\xi}_d$ and degrees $\tilde{p}_d$ with $d=1,\hdots,D$. The weight function $\omega(\hat{x})$ of the mass matrix is interpolated as
\begin{equation}
\omega(\hat{x}) \approx \sum_{\mathbf{j}\in \mathcal{J}} W_{\mathbf{j}}\tilde{\beta}_{\mathbf{j}}(\hat{x})=W:\tilde{B}(\hat{x}),
\end{equation}
where $\tilde{\beta}_\mathbf{j}(\hat{x})$ are the elements of the spline space $\tilde{\mathbb{S}}^D$ and $\tilde{B}(\hat{x})$ is the tensor holding all $\beta_\mathbf{j}$ ordered according to the index set $\mathcal{J}$.
The weight tensor $W$ has the same dimension as the spline space $\tilde{\mathbb{S}}^D$, being $(\tilde{n}_1,\hdots,\tilde{n}_D)$, and we get its entries by interpolating the weight function in a sufficient number of points, namely $\tilde{n} = \tilde{n}_1\cdots\tilde{n}_D$.

We construct canonical low rank representations of the weight tensor,
\begin{equation} \label{equation:canonicalLowRank}
W \approx \sum_{r=1}^R \bigotimes_{d=1}^D w_r\hk =: W_R,
\end{equation}
with $w_r\hk \in \R^{n_d}$ to get a low rank representation of the weight function,
\begin{equation} \label{equation:weightLowRank}
\omega(\hat{x}) \approx W_R : \tilde{B}(\hat{x}) = \sum_{r=1}^R \prod_{d=1}^D w_r\hk \cdot \tilde{B}\hk(\hat{x}\hk).
\end{equation}
Here $\tilde{B}\hk(\hat{x}\hk) \in \R^{n_d}$ denotes the vector holding all univariate basis functions evaluated in $\hat{x}\hk$ as in \cref{eq:Bd}, and ``$\cdot$'' is the scalar product.
The entries of the mass matrix can be approximated using this low rank representation and we can calculate each entry as the sum of products of univariate integrals,
\begin{align}
M_{\mathbf{i},\mathbf{j}} = \sum_{r=1}^R \prod_{d=1}^D  \int_0^1\beta_{i_d}\hk\beta_{j_d}\hk w_r\hk \cdot \tilde{\beta}\hk \D \hat{x}\hk.
\end{align}
With these univariate integrals we define a univariate mass matrix, which depends on some weight function $\omega$, as
\begin{equation} \label{equation:massMatrix1D}
M\hk(\omega) = \int_0^1 B\hk \otimes B\hk \omega \, \D \hat{x}\hk,
\end{equation}
where $B\hk\in \R^{n_d}$ is the vector holding all $n_d$ univariate B-splines of $\mathbb{S}^{p_d}_{\xi_d}$.
According to the tensor representation in \Cref{equation:massTensor1}, we can finally write the mass matrix as a sum of Kronecker products of small univariate mass matrices \cref{equation:massMatrix1D} with $\omega = w_r\hk \cdot \tilde{\beta}\hk$,
\begin{equation} \label{equation:massFinal}
M = \sum_{r=1}^R \bigotimes_{d=1}^D M\hk(w_r\hk \cdot \tilde{\beta}\hk).
\end{equation}

The same procedure can be applied to the weight function of the stiffness matrix $Q(\hat{x})$. Note that $Q(\hat{x}) \in \R^{D\times D}$, thus we have to apply the interpolation to each entry of $Q$. Similarly to \cref{equation:weightLowRank}, for each entry of $Q$ we get the canonical low rank representation
\begin{equation}
q_{k,l}(\hat{x}) \approx V_{k,l,R} : \tilde{B}(\hat{x}) = \sum_{r=1}^R \prod_{d=1}^D v_{k,l,r}\hk \cdot \tilde{\beta}\hk (\hat{x}\hk), \quad \mbox{ for all } k,l=1,\hdots,D,
\end{equation}
with $v_{k,l,r}\hk \in \R^{n_d}$.

Using this low rank method, we approximate the entries of the stiffness matrix as
\begin{align} 
K_{\mathbf{i},\mathbf{j}} &= \sum_{k,l=1}^D \int_{\hat{\Omega}} \Big ( \prod_{d=1}^D \delta(l,d) \beta_{i_d}\hk \delta(k,d) \beta_{j_d}\hk \Big )\sum_{r=1}^R \prod_{d=1}^D v_{k,l,r}\hk \cdot \tilde{\beta}\hk \D \hat{x},\\
& = \sum_{k,l=1}^D \sum_{r=1}^R \prod_{d=1}^D \int_0^1 \delta(l,d) \beta_{i_d}\hk \delta(k,d)\beta_{j_d}\hk v_{k,l,r}\cdot \hk\tilde{\beta}\hk \D \hat{x}\hk,
\end{align}
where $\mathbf{j}=(j_1,\ldots,j_D)$, and $\delta(k,d)$ denotes the operator acting on $f$ as
\begin{equation}
\delta(k,d) f = \begin{cases} \frac{\partial f}{\partial \hat{x}_d} &\mbox{ if } k = d, \\ f & \mbox{ otherwise}. \end{cases}
\end{equation}
To get a representation for the stiffness matrix corresponding to the mass matrix representation in \cref{equation:massFinal}, we define the $D^2$ univariate stiffness matrices dependent on some weight function $q^{(d)}(\hat{x}\hk)$ as
\begin{equation}
K_{k,l}\hk(q^{(d)}) = \int_0^1 \left (\delta(l,d) B \right ) \otimes \left (\delta(k,d) B \right ) q^{(d)} \D \hat{x}\hk, \quad \mbox{ for } k,l=1,\hdots,D.
\end{equation}
With this and $q^{(d)} = v_{k,l,r}\hk\cdot \tilde{\beta}\hk$ the final low rank tensor representation of the stiffness matrix is
\begin{equation} \label{equation:stiffnessFinal}
K = \sum_{k,l=1}^D \sum_{r=1}^R \bigotimes_{d=1}^D K_{k,l}\hk(v_{k,l,r}\hk \cdot \tilde{\beta}\hk).
\end{equation}

Strategies to efficiently compute the low rank representation for $W$ have been discussed in \cite{angelos2} for two dimensional settings and for three dimensional settings with partial decompositions in \cite{angelos3} and using tensor decompositions in \cite{angelos3,buenger2020}.

We follow the strategy presented in \cite{buenger2020} using low rank Tensor Train methods \cite{osel-tt-2011,dolgov:savostyanov:2014} for the decompositions and further computations as it allows us to exploit the resulting low-rank structure.

 \section{The Alternating Minimal Energy solver for Tensor Train} \label{chapter:amen}
 The IAS algorithm requires solving the minimization problem \cref{eq:LS} multiple times in line 4 of \cref{alg:local,alg:global}. To solve the KKT system \eqref{eq:saddlepoint} we use a Block-structured Alternating Minimal Energy solver (AMEn) as proposed in \cite{buenger2020}. The saddle point system is large and each block has a Kronecker product based structure. This structure can be interpreted as a low-rank Tensor-Train (TT) \cite{osel-tt-2011} representation. Alternatively, using a CGLS method with the reduced formulation \eqref{eq:reduced} requires a large number of matrix-vector products with the matrix $A$. This matrix is dense and quickly gets unfeasibly large. Therefore, we want to avoid its complete computation and instead solve the matrix-vector product without explicitly forming $A$. We propose using the TT AMEn method to apply this computation as well.
 
 The TT format is especially convenient for our purpose as the underlying Kronecker product structure of the parameter-to-observable map \eqref{eq:kronecker1} -- \eqref{eq:kronecker4} can be interpreted as a low-rank Tensor-Train representation. Using the TT format we can easily compute the desired result in a low-rank format without exceeding the memory limitations of a standard computer even for large models, which usually are infeasible.
 
 The Tensor-Train format represents a $d$-dimensional tensor $\mathbf{T} \in \mathbb{W}_{n_1,n_2,\hdots, n_d}$ of order $(n_1,n_2,\hdots,n_d)$ with so-called TT-cores $T^{(k)} \in \R^{r_{k-1} \times n_k \times r_k}$ with TT-ranks $r_k$ for $k = 1, \hdots,d$. By convention we set $r_0 = r_d =1$. Each TT-core can be interpreted as a parameter-dependent matrix $T^{(i)}(j_k) \in \R^{r_{k-1} \times r_k}$, $j_k = 1,\hdots, n_i$ and every element of $\mathbf{T}$ is represented as the product
 \begin{equation}
 \mathbf{T}(j_1,\hdots,j_d) = T^{(1)}(j_1) \cdots T^{(d)}(j_d).
 \end{equation}
 
 The whole tensor can be written as a sum of Kronecker products,
 \begin{equation}
 \mathbf{T} = \sum_{\alpha_1 = 1}^{r_1} \cdots \sum_{\alpha_d = 1}^{r_d} \bigotimes_{k=1}^{d} \mathbf{T}_{\alpha_{k-1}, \alpha_k}^{(k)} =: tt  (T^{(1)},\hdots,T^{(d)}  ) ,
 \end{equation}
 where the subscripts $\alpha_{k-1}$, $\alpha_k$ are row and column indices of the TT-core $T^{(k)}$. This format directly corresponds to the matrices in \eqref{eq:kronecker1} -- \eqref{eq:kronecker4}.
 
 We want to solve an equation system $\mathbf{A} \mathbf{x} = b$ in this format, where $\mathbf{A}$ is a tensor of size $N \times N$ with $N = (n_1, \hdots, n_d)$ and $\mathbf{x}$ and $b$ are tensors of size $N$.
 This can be done efficiently with an energy function minimization cycling over the TT-cores, the \textit{alternating linear scheme} (ALS) \cite{HoltzRohwedderSchneider:2012} . This approach constructs low-dimensional systems of linear equations for each core, which can then be solved with standard numerical methods. 
 Here, we derive the method for symmetric $\mathbf{A}$. Solving $\mathbf{A} \mathbf{x} = b$ then corresponds to the minimization of the energy function,
 \begin{equation}
 \min_{\mathbf{x}} J(\mathbf{x}) = \|\mathbf{x}_* - \mathbf{x}\|^2_\mathbf{A} =  (\mathbf{x}, \mathbf{A}\mathbf{x}) - 2 (\mathbf{x}, \mathbf{b}) + \text{const}, \label{eq:energy}
 \end{equation}
 with the exact solution $\mathbf{x}_* = \mathbf{A}^{-1}\mathbf{b}$. To find a solution for \eqref{eq:energy} we make an initial guess $\mathbf{x}_0$ and cycle over its TT-cores where we solve a local problem to improve the current guess. For this, in iteration $k$ all cores but the $k$-th are \textit{frozen} and we minimize over
 \begin{align}
 \mathbf{x}_\text{new} &= tt\left (x^{(1)},\hdots, x^{(k-1)},x_\text{new}^{(k)},x^{(k+1)},\hdots,x^{(d)} \right ) \\
 x_\text{new}^{(k)} &= \min_{x^{(k)}} J(\mathbf{x}_\text{new}).
 \end{align}
 The energy function $J(\mathbf{x})$ does not grow during updates and the solution will subsequently converge to a local minimum. The tensor-train format is linear in its cores, as
 \begin{equation}
 \mathbf{x} = tt \left ( x^{(1)}, \hdots, x^{(k)},\hdots, x^{(d)} \right ) = \mathbf{x}_{\neq k} x_k,
 \end{equation}
 where $\mathbf{x}_{\neq k}$ is the Tensor-Train where the $k$-th core $\mathbf{x}^{(k)}  \in \R{r_{k-1}\times n_k \times r_k}$ is replaced by an identity operator of the same size and $x_k$ is a vectorization of $x^{(k)}$. With this, the energy function for the local problem becomes
 \begin{equation}
 J(\mathbf{x}) = (\mathbf{A}\mathbf{x}_{\neq k}x_k, \mathbf{x}_{\neq k}x_k) - 2(\mathbf{b},\mathbf{x}_{\neq k}x_k) = (\mathbf{x}_{\neq k}^*\mathbf{A}\mathbf{x}_{\neq k}x_k,x_k) - 2(\mathbf{x}_{\neq k}^*\mathbf{b},x_k). \label{eq:local:energy}
 \end{equation}
 The gradient of \eqref{eq:local:energy}  with respect to $x_k$ is zero when
 \begin{equation}
 (\mathbf{x}_{\neq k}^*\mathbf{A}\mathbf{x}_{\neq k})x_k = \mathbf{x}_{\neq k}^* y. \label{eq:local:system}
 \end{equation}
 Therefore, the solution to the local minimization is equal to a solution in a reduced basis. This problem is small and can be solved by standard numerical methods. For the general ALS approach the resulting TT-ranks - and therefore the maximum accuracy - are fixed by the ranks set in the initial guess. But with some extension we can adapt the TT-ranks of the solution dynamically. For the purpose of calculating a low-rank approximation, we chose the so-called \textit{Alternating Minimal Energy} (AMEn) method \cite{dolgov:savostyanov:2014} . 
 \begin{algorithm}
 	\caption{Low-rank AMEn} \label{alg:amen}
 	\begin{algorithmic}[1]
 		\State Input: Initial guess $\mathbf{x}$ with TT-ranks $r = (r_1,\hdots,r_d)$, tolerance $\epsilon_0$ and/or rank bounds $r_\text{max}$, system in TT-format $\mathbf{A}$, $\mathbf{b}$
 		\While{$\| \mathbf{A} \mathbf{x} - \mathbf{b} \| > \epsilon_0$ and $k < \text{maxIter}$}
 		\State Form $A_1 = \mathbf{x}_{\neq 1}^* \mathbf{A} x_{\neq 1}$, $b_1 = \mathbf{x}_{\neq 1}^* \mathbf{b}$
 		\State Solve $A_1 u^{(1)} = b_1$
 		\State Optional: Do an SVD-based compression of $u^{(1)}$
 		\State Let $\mathbf{u} = tt(u^{(1)},x^{(2)},\hdots, x^{(d)})$ and $\mathbf{r} =  \mathbf{A} \mathbf{u} - \mathbf{b}$
 		\If{$\|r\| < \epsilon_0$ }
 		\State \Return $\mathbf{x} = \mathbf{u}$
 		\EndIf
 		\State Find smallest-rank residual approximation $\mathbf{\bar{r}}$ with $\|\mathbf{\bar{r}} -\mathbf{r}\| \leq \epsilon_0 \|\mathbf{r}\|$ and/or $\text{rank}(\bar{r}^{(1)}) \leq r_\text{max}$
 		\State Expand basis by first TT-core $\bar{r}^{(1)}$ of $\mathbf{\bar{r}}$, $x^{(1)} = \begin{bmatrix} u^{(1)} & \bar{r}^{(1)} \end{bmatrix}$, and accordingly $x^{(2)} = \begin{bmatrix} x^{(2)} \\ 0 \end{bmatrix}$
 		\State Form $X_1 = x^{(1)} \otimes I_{n_2,\hdots,n_d}$ of size $(n_1,\cdots,n_d)\times (r_1n_2,\cdots,n_d)$
 		\State Form the $(d-1)$-dimensional system $A_{\geq 2} x_{\geq 2} = b_{\geq 2}$ with $A_{\geq 2} = X_1^* A X_1$, $b_{\geq 2} = X_1^* b$.
 		\If{d = 2}
 		\State Solve $A_{\geq 2} x_{\geq 2} = b_{\geq 2}$ directly
 		\Else
 		\State Call AMEn recursively with $\mathbf{x} \leftarrow tt(x^{(2)},\hdots,x^{(d)})$, $\mathbf{A} \leftarrow A_{\geq 2}$, $\mathbf{b} \leftarrow b_{\geq 2}$.
 		\EndIf
 		\EndWhile
 		\State Output: low-rank solution $\mathbf{x} = tt(x^{(1)},\hdots,x^{(d)})$
 	\end{algorithmic}
 \end{algorithm}
 
 The algorithm is outlined in \cref{alg:amen}. This method first solves a local system as in \cref{eq:local:system}. Then it expands the components of the solution subsequently by local gradient information as in line 11 of the algorithm. Afterwards, the system dimension is reduced by one and the previous computation repeated in a recursive way until a 2D matrix equation system is reached. Solving this equation system finalizes the outer iteration after which all steps are repeated until reaching the convergence tolerance. The method proved to be robust and has a fast convergence rate. For a detailed analysis and more information we refer the reader to \cite{dolgov:savostyanov:2014}.
 
 Using the TT notation and the AMEn method, we can compute matrix-vector products with $A$ efficiently without forming the full matrix and get the resulting solution vectors in a low-rank format. For solving the block system \cref{eq:saddlepoint} we use a block solver based on AMEn \textit{AMEn Block Solve}, which exploits the block structure of \cref{eq:saddlepoint} as well as the low rank tensor product structure of each component.
 
 \section{Numerical experiments} \label{chapter:examples}
We want to illustrate the performance of our method with some numerical examples and analyze the choice of parameters for PDE problems. We first present results for different parameter settings for a two dimensional domain and later show results for a 3D model. The parameters for the different hypermodel setups are listed in \cref{table:parameters}.

\begin{table}
	\centering
	\begin{tabular}{l||r|r|r}
		$r$ & $\eta$ & $\beta$ & $\vartheta$ \\
		\hline
		1 & $10^{-5}$ & $3/2 + \eta$ & 3.3 \\
		0.5 & $10^{-3}$ & $(3/2 + \eta)/r$ & 8.3 \\
		$-1$ & $\beta + 3/2$ & 3 & $1.5 \times 10^{-4}$
	\end{tabular}
	\caption{Parameters for the three hyper models} \label{table:parameters}
\end{table}

\subsection*{Example 1: 2D domain}
\begin{figure}[h]
	\centering
	\begin{subfigure}[b]{0.47\textwidth}
	\includegraphics[width=1\textwidth]{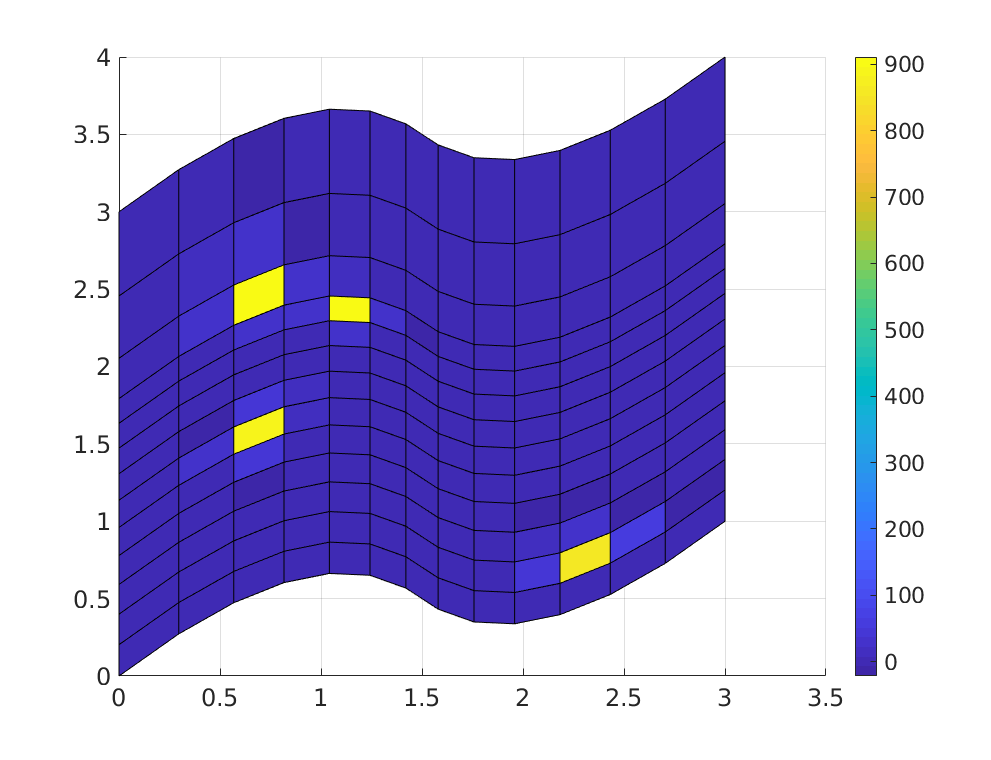}
	\caption{Two dimensional domain with sparse heat sources.} \label{fig:2D}
	\end{subfigure}\hfill%
	\begin{subfigure}[b]{0.47\textwidth}
		\includegraphics[width=0.9\textwidth]{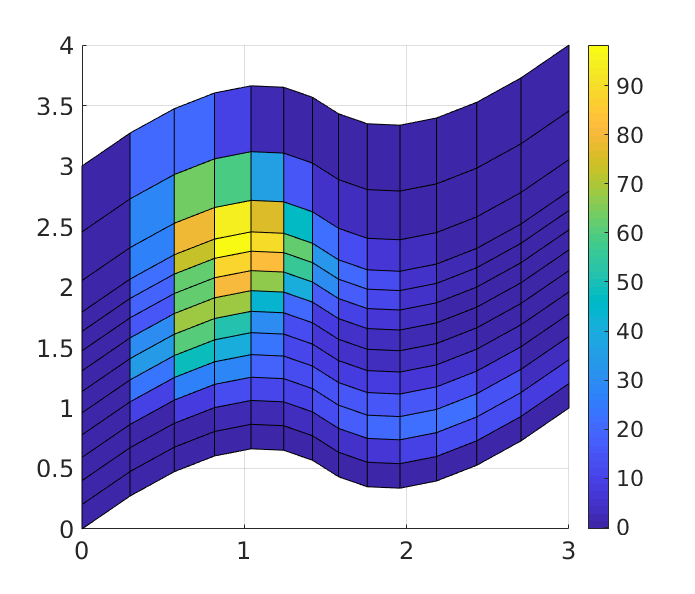}
		\caption{Two dimensional domain after 100 time steps at end time $T = 1$ with noise.} \label{fig:endState}
		\end{subfigure}
\begin{subfigure}[t]{0.47\textwidth}
	\includegraphics[width=0.9\textwidth]{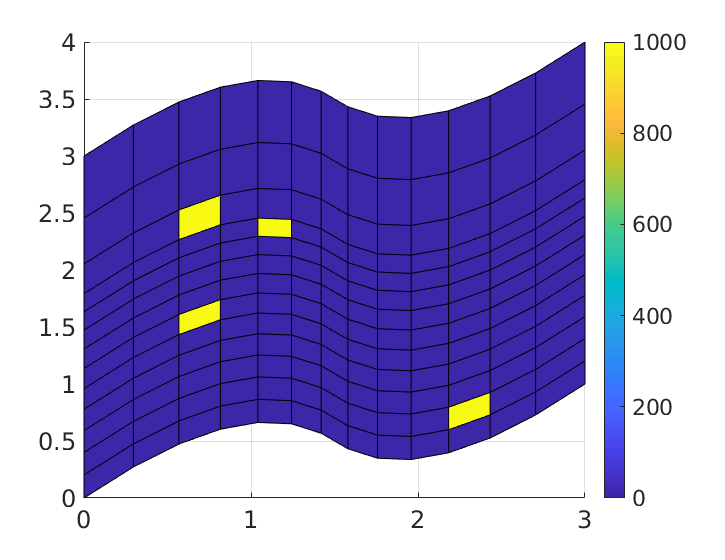}
	\caption{Resulting sparse reconstruction.} \label{fig:result}
\end{subfigure} \hfill%
\begin{subfigure}[t]{0.47\textwidth}
	\includegraphics[width=0.9\textwidth]{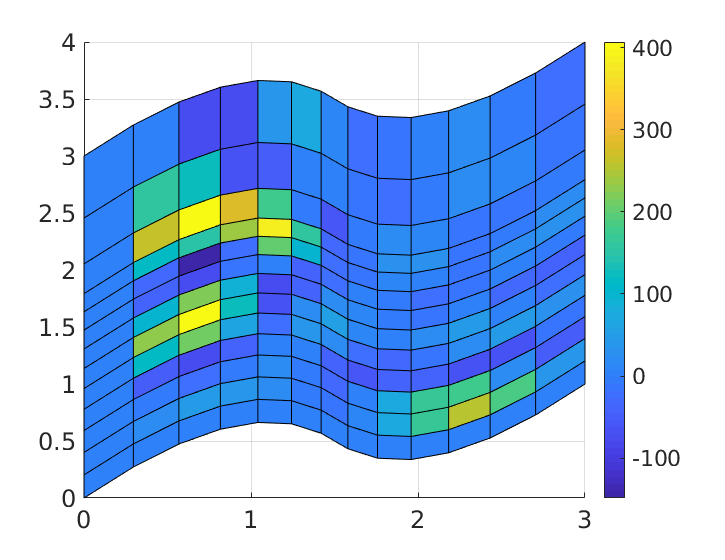}
	\caption{For comparison: Resulting reconstruction using only least squares without the hyperprior pdf.} \label{fig:LSresult}
\end{subfigure}
\caption{Setup for the 2D example problem. The PDE is discretized using 1024 basis functions.}
\end{figure}

First, let us consider a two dimensional domain with sparse heat sources, as shown in \cref{fig:2D}. This domain is discretized using IGA with 32 B-splines per dimension resulting in 1024 spatial degrees of freedom. The total rank for the mass matrix is 8 and for the stiffness matrix we have a combined number of 12 low-rank components. For the time discretization we use $N_t = 50$.

This setup is governed by a heat distribution 
\begin{align}
\frac{\partial x}{\partial t}(t) &= 0.1 \Delta x(t)  &&\mbox{ in } \Omega \times (0,1] \\
x(0) & = x_0  &&\mbox{ in } \Omega \times (t=0)\\
x(t) &= 0 && \mbox{ on } \partial \Omega \times [0,1]
\end{align} 
and we measure the state only once after 50 time steps at $T = 1$. The resulting state is shown in \cref{fig:endState} and is measured under noise $e \sim \mathcal{N}(0,0.1)$. We run a non-modified IAS algorithm with hyperprior parameters as in the first row of \cref{table:parameters} with $r = 1$, as suggested in \cite{calvetti2020}. We do 50 IAS steps and 30 CGLS steps per iteration. For the low-rank matrix vector products we set a tolerance of $10^{-6}$. The solver tolerance for AMEn was set to $10^{-4}$ and we do a maximum of 20 AMEn steps per iteration. An exemplary result is shown in \cref{fig:result}. We see that the solution indeed fits the very sparse initial state perfectly. The algorithm successfully reconstructs the 4 sources. For comparison the result for a least squares reconstruction is shown in \cref{fig:LSresult}. 

\subsection*{Example 2: Hyperprior parameters}
Let us now take a closer look at the parameters of the generalized gamma distribution used for the hyperprior in \Cref{eq:gamma} and the performance of the local and global hybrid IAS algorithm. We use the same setup as before.
We set up a globally convex hyperprior model with $r = 1$ two different sparsity inducing models with $r = 0.5$ and $r = -1$ to compare the performance of the different algorithm variations. The rest of the parameters correspond to \cref{table:parameters}.

We will compare three different methods: CGLS with a full computation of $A$ denoted by full CGLS that will not be feasible for large scale problems; CGLS with low-rank matrix-vector products in TT format denoted by TT CGLS; and solving the KKT system of the constrained optimization problem with the block AMEn method denoted by AMEn.

\begin{figure}[h]
	\centering
	\begin{subfigure}[t]{0.49\textwidth}
		\includegraphics[width=1\textwidth]{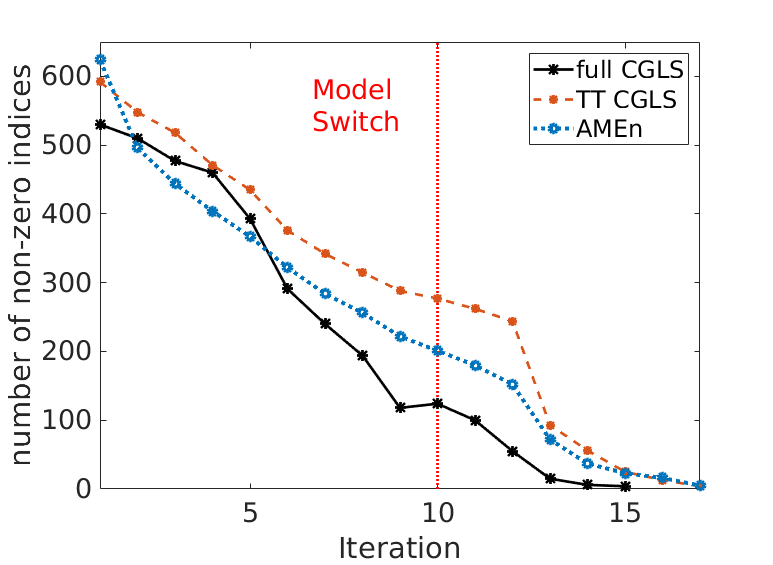}
		\caption{Number of indices for which $u_0$ is non-zero per iteration for $r=-1$. The switch to the model with stronger sparsity promotion is indicated by the red dotted line.} \label{fig:outputglobal-1}
	\end{subfigure}\hfill
	\begin{subfigure}[t]{0.49\textwidth}
		\includegraphics[width=1\textwidth]{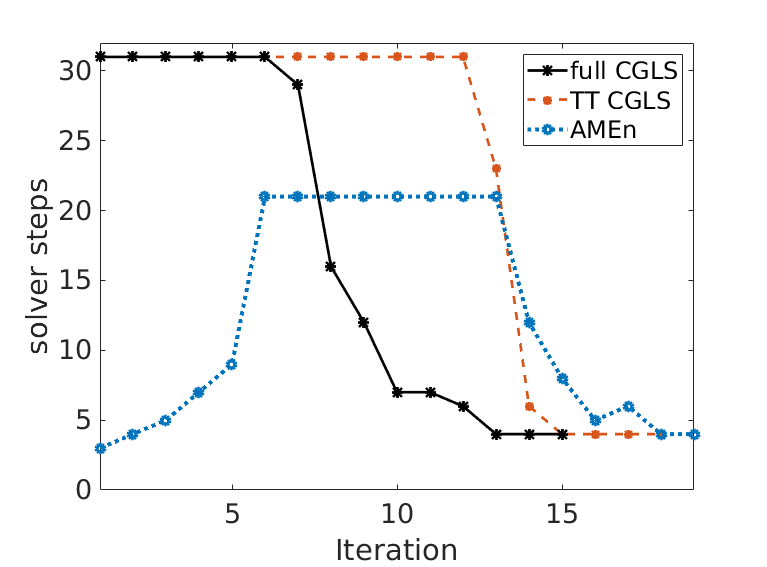}
		\caption{CGLS steps and AMEn sweeps per iteration for the global hybrid IAS with $r=-1$.} \label{fig:CGLSglobal-1}
	\end{subfigure}
	\begin{subfigure}[t]{0.49\textwidth}
		\includegraphics[width=1\textwidth]{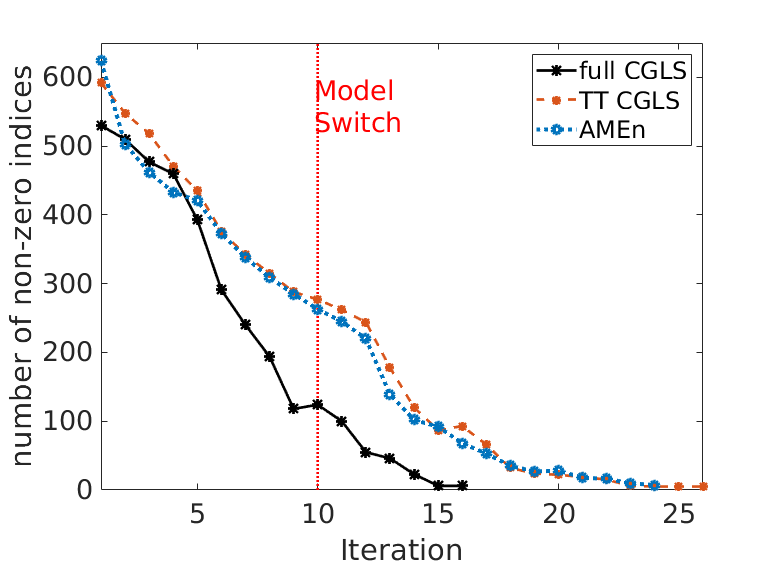}
		\caption{Number of indices for which $u_0$ is non-zero per iteration for $r=0.5$. The switch to the model with stronger sparsity promotion is indicated by the red dotted line.} \label{fig:outputglobal05}
	\end{subfigure}\hfill
	\begin{subfigure}[t]{0.49\textwidth}
		\includegraphics[width=1\textwidth]{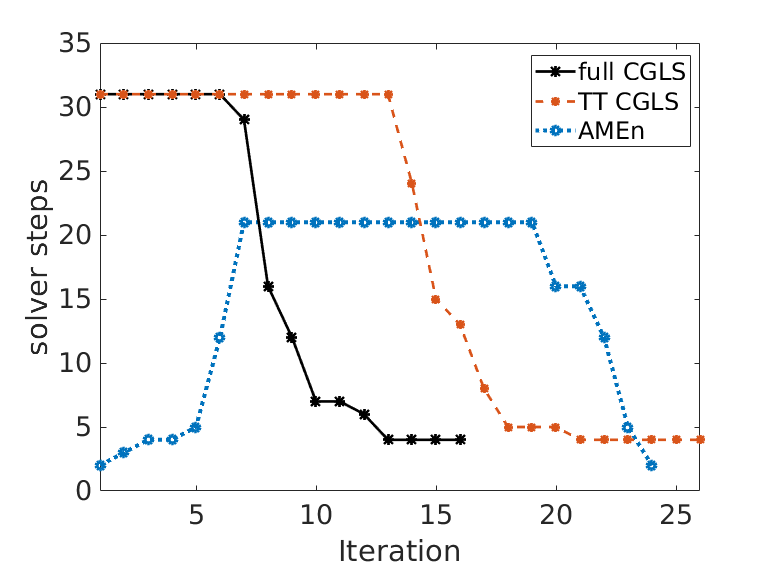}
		\caption{CGLS steps and AMEn sweeps per iteration for the global hybrid IAS with $r=0.5$.} \label{fig:CGLSglobal05}
	\end{subfigure}
	\caption{Development of non-zero indices and CGLS steps for the global hybrid IAS. We compare using the full matrix $A$ with CGLS (full CGLS), CGLS with low-rank matrix-vector products (TT CGLS) and a constrained minimization with AMEn.} \label{fig:globalIAS}
\end{figure}

First, let us review the global hybrid IAS switching models after 10 iterations. In \cref{fig:globalIAS} we see results for the two different hypermodel setups with $r=-1$ and $r=0.5$. We compare the performance of the low-rank CGLS and the constrained optimization with AMEn. Both setups start out with a globally convex hypermodel with $r = 1$. After 10 iterations the first algorithm switches to a greedy hypermodel with $r=-1$, indicated by a red vertical line. We can see the number of non-zero entries reducing shortly after the tenth iteration in \cref{fig:outputglobal-1} for $r=-1$. In \cref{fig:CGLSglobal-1}, we see that the CGLS method starts with the maximum number of iterations and the low-rank version requires more IAS steps until the number of local iterations starts reducing compared to the full CGLS method. This can be attributed to the inexact matrix-vector products leading to inexact descent directions. Nonetheless, the low-rank CGLS converges with only a small number of extra iterations required. 

The constrained optimization method does not require the full number of 20 steps from the beginning. This is due to $D_{\theta}^{-1/2}$ starting out quite large. This could be interpreted as a large regularization parameter for $u_0$ in the optimization, which is generally easier to solve. As $\theta$ changes, the minimization gets more challenging before the iteration numbers start decreasing again. All three algorithms converge with the local solvers needing 4 steps for their last iterations, which corresponds to the number of non-zero values in $u_0$.

In \cref{fig:outputglobal05,fig:CGLSglobal05} we see very similar results for a second model with $r=0.5$. Note, that both algorithms are equal until iteration 10 where we switch to the different second models. For $r=0.5$ convergence is not reached as fast as for $r=-1$, which was to be expected as $r=-1$ promotes sparsity more strongly. Again we see a quick reduction in non-zero entries shortly after the switch to the second hyper model in \cref{fig:outputglobal05}. And as before the low-rank methods need some more iterations to converge, as seen in \cref{fig:CGLSglobal05}. 

\begin{figure}[h]
	\centering
	\begin{subfigure}[t]{0.49\textwidth}
		\includegraphics[width=1\textwidth]{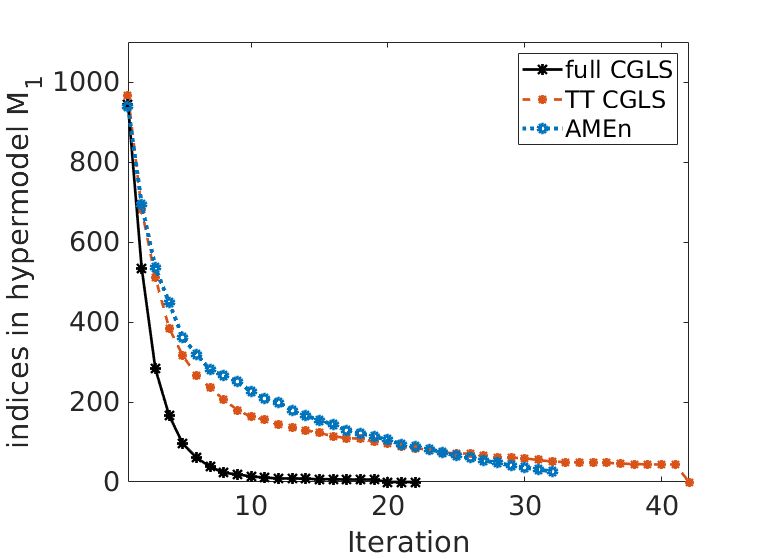}
		\caption{Number of indices remaining in hypermodel $M_1$ in each iteration for the local hybrid IAS with inverse gamma $r=-1$. } \label{fig:outputlocal-1}
	\end{subfigure}\hfill
	\begin{subfigure}[t]{0.49\textwidth}
		\includegraphics[width=1\textwidth]{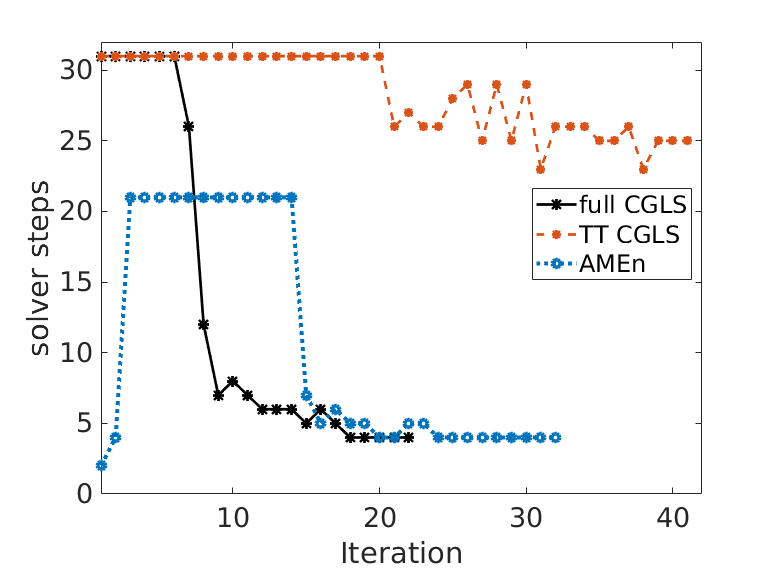}
		\caption{CGLS and AMEn steps per iteration for the local hybrid IAS with $r=-1$.} \label{fig:CGLSlocal-1}
	\end{subfigure}
	\begin{subfigure}[t]{0.49\textwidth}
		\includegraphics[width=1\textwidth]{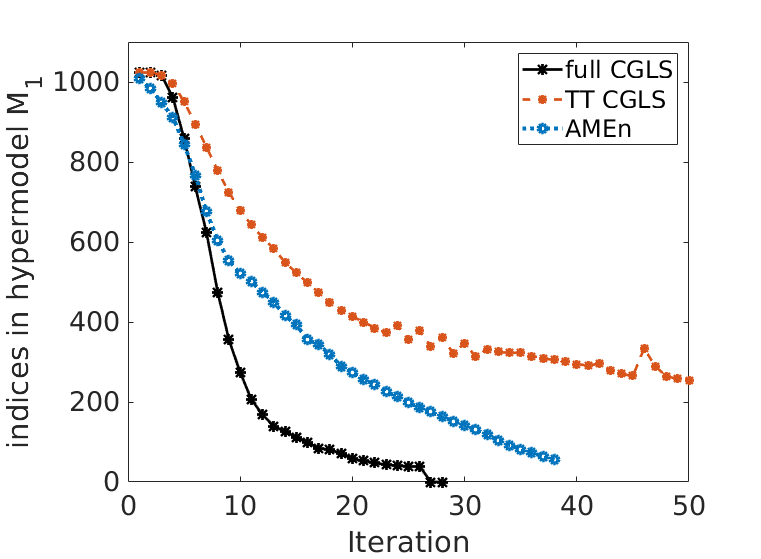}
		\caption{Number of indices remaining in hypermodel $M_1$ in each iteration for the local hybrid IAS with less greedy gamma parameter $r=0.5$.} \label{fig:outputlocal05}
	\end{subfigure}\hfill
	\begin{subfigure}[t]{0.49\textwidth}
		\includegraphics[width=1\textwidth]{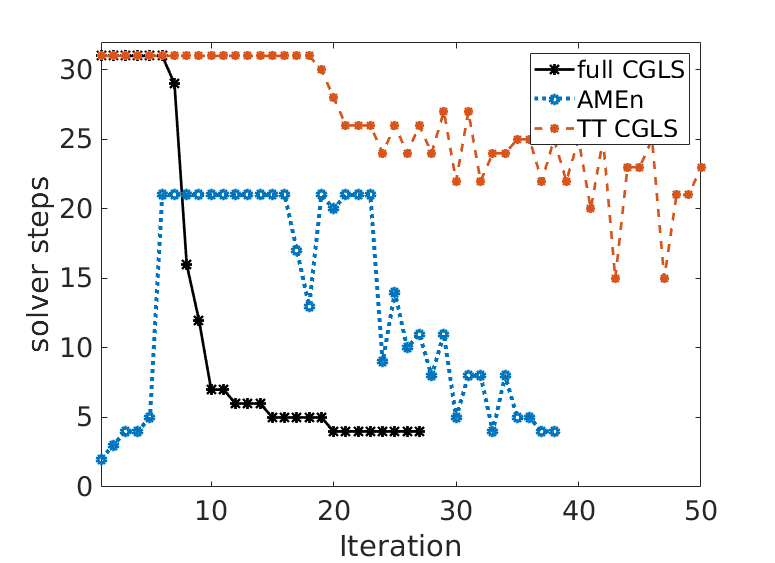}
		\caption{CGLS and AMEn steps per iteration for the local hybrid IAS with $r=0.5$. CGLS with TT did not reach convergence.} \label{fig:CGLSlocal05}
	\end{subfigure}
	\caption{We compare using the full matrix $A$ with CGLS (full CGLS), CGLS with low-rank matrix-vector products (TT CGLS) and a constrained minimization with AMEn.} \label{fig:localIAS}
\end{figure}
Next, we look at the local hybrid IAS. Here, the switch to the second models is done individually for each entry of $\theta$ as described in \cref{chapter:hyperpriors}. The results are shown in \cref{fig:localIAS}. The first row shows the performance with $r=-1$ for the second model. In \cref{fig:outputlocal-1} we see the number of entries in $\theta$, which remain in the first hypermodel with $r=1$. All others are switched to $r=-1$ as soon as their convexity condition is fulfilled. We see that the constrained minimization with AMEn performs quite well. The number of parameters remaining in the model with $r=1$ is reducing quickly. The low-rank CGLS requires some more iterations but ultimately converges as well. In \cref{fig:CGLSlocal-1}, we see that again the AMEn method ends up converging with only 4 local iterations towards the end of the algorithm, similar to the full CGLS method. The low-rank TT CGLS however requires a large number of steps until the very end when it converges.

 Additionally, in \cref{fig:outputlocal05,fig:CGLSlocal05} we see the results for a model with $r=0.5$. As expected the method needs more iterations to converge and the number of indices in $\theta$ remaining in the first model decreases more slowly for all three algorithms in \cref{fig:outputlocal05}. Again, the full CGLS and AMEn converge with 4 local iterations towards the last iterations of the algorithms. However, the TT CGLS does not converge within 50 outer iterations. Here, a smaller residual tolerance for the matrix-vector products might result in more successful steps.
 
 \begin{table}
 	\centering
 	\begin{tabular}{l l || r r| r r  }
 		& & \multicolumn{2}{c}{global} & \multicolumn{2}{c}{local}  \\
 		& & iterations & residual & iterations & residual \\
 		\hline
 		full CGLS & $r=-1$ & 15 &2.3e-7 & 22 & 9.8e-3  \\
	 		& $r=0.5$ & 16 & 5.9e-6 &28 & 7.2e-3 \\
 		TT CGLS & $r=-1$ & 18 & 9.2e-7 & 41 & 4.8e-2 \\
	 		& $r=0.5$ & 26 & 2.6e-6 & $>50$ & 0.4 \\
 		AMEn & $r=-1$ & 18 & 3.6e-4 & 32 & 8.3e-2 \\
	 		& $r=0.5$ & 24 & 2.4e-4 & 38 & 8.2e-2
 	\end{tabular}
 	\caption{Number of total iterations and residual error to the original initial state as $\|u_0 - u\|/\|u_0\|$ with true initial state $u_0$ and reconstruction $u$.} \label{table:2D}
 \end{table}

The number of iterations and the resulting residual with respect to the actual initial state are listed in \cref{table:2D}. We see that overall the results for the global method are very good. The larger residuals for AMEn correspond to the solver tolerance, which was set at $10^{-4}$. For the local method however all methods required a larger number of iterations and the resulting reconstruction was not as close to the original as for the global method. 

\subsection*{Example 2: 3D domain}

\begin{figure}[h]
	\centering
	\begin{subfigure}[b]{0.95\textwidth}
		\includegraphics[width=0.95\textwidth]{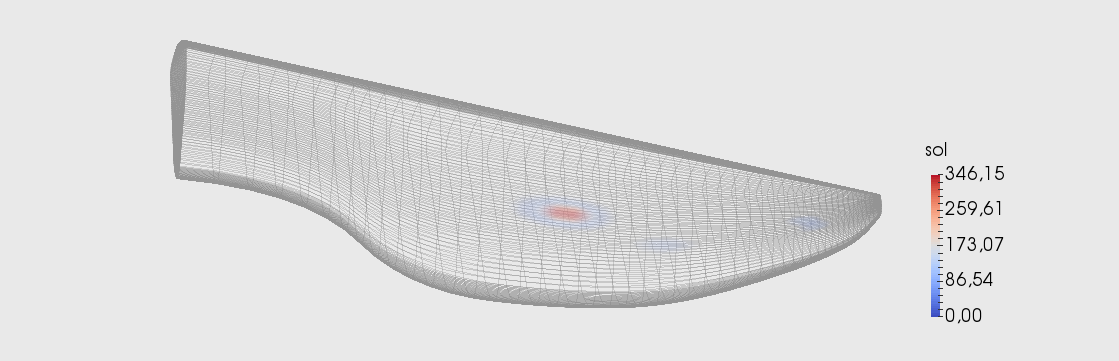}
		\caption{Three dimensional domain with sparse heat sources.} \label{fig:3D}
	\end{subfigure}
	\hfill
	\begin{subfigure}[b]{0.95\textwidth}
		\includegraphics[width=0.95\textwidth]{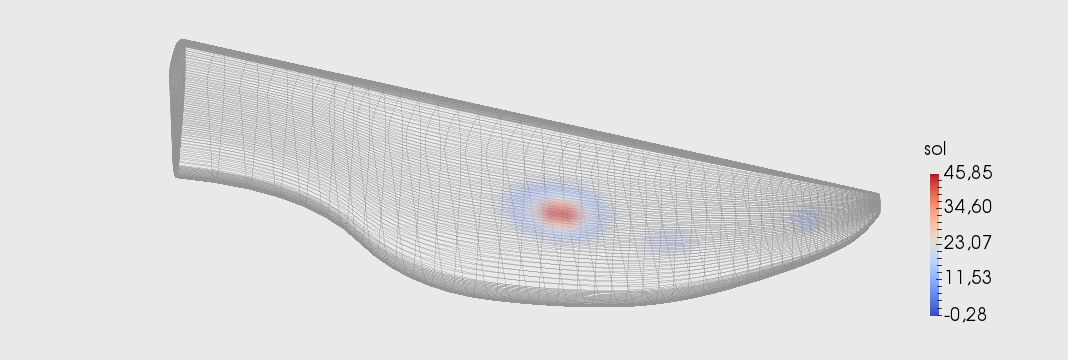}
		\caption{Three dimensional domain after 100 time steps at end time $T = 1$ with noise.} \label{fig:endState3D}
	\end{subfigure}
	\caption{Setup for the 3D example problem} \label{fig:3Dproblem}
\end{figure}

\begin{figure}
	\centering
		\includegraphics[width=0.95\textwidth]{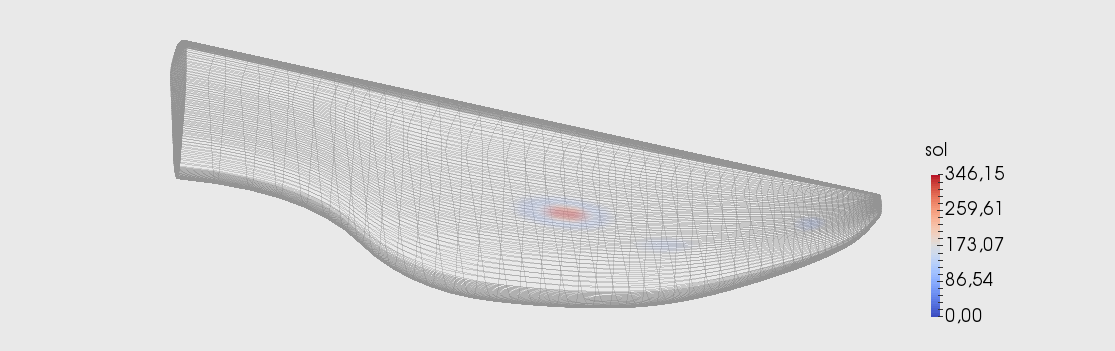}
	\caption{Resulting sparse reconstruction for the 3D example problem}\label{fig:result3D}
\end{figure}

We now consider a three dimensional domain with a more complex geometric model. The low-rank structure has a much higher rank than the previous example, namely 10 factors for the mass and 180 factors for the stiffness matrix. The initial state is displayed in \cref{fig:3D} and the final state we use for the reconstruction with noise $\sigma = 0.1$ in \cref{fig:endState3D}. An exemplary reconstruction with global hybrid IAS and $r_2 = -1$ is displayed in \cref{fig:result3D}.

We modify the problem size by changing the space discretization size and number of time steps to observe the scaling behavior of the method. For the TT CGLS method we set the tolerance to 1e-6. The result for two spatial discretizations and three different time resolutions is displayed in \cref{figure:3Dsteps}. The setups with equal spatial discretization behave very similarly. For a total of 756 degrees of freedom all time discretizations result in less than 10 iterations. For the discretization with 6048 spatial nodes some more iterations are required but again the method displays similar behavior for all time discretizations. For the constrained minimization with AMEn the behavior is quite different. Here, we display three spatial discretizations and two time discretizations in \cref{table:3Dsweeps}. For the different spatial resolutions the number of IAS iterations is somewhat stable. But with a larger number of time steps the required iterations increase drastically. In future research, we aim to mitigate this behavior by further tailoring the AMEn solver regarding the enrichment process and the design of suitable preconditioners.

\begin{table}[h]
	\centering
	\begin{minipage}{0.58\textwidth}
	\includegraphics[width=0.99\textwidth]{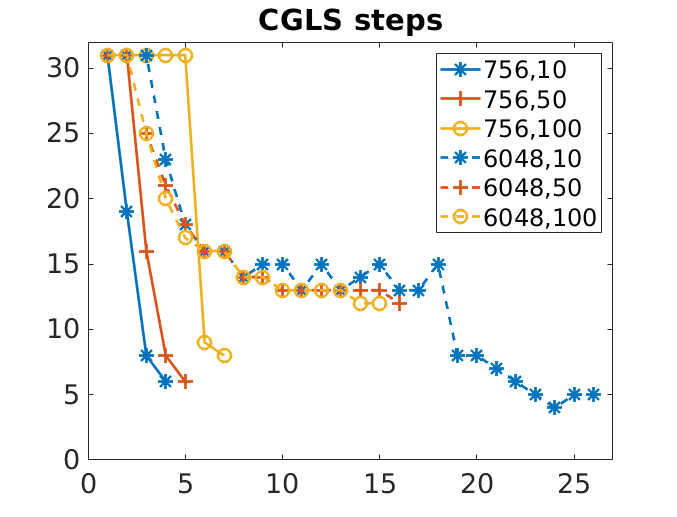}
	\subcaption{Number of IAS steps for the global hybrid IAS with $r=-1$ and the 3D model problem for different discretization sizes. Here we used TT CGLS} \label{figure:3Dsteps}
\end{minipage}\hfill%
\begin{minipage}{0.37\textwidth}
	\centering
	\begin{tabular}{l||r|r|r}
		$N_t \backslash n$ & 756 & 6048 &  48384\\
		\hline 
		50 & 16 & 18 & 29\\
		100 & 47 & 50 & 49 \\
	\end{tabular}
	\caption{Number of IAS steps using a constrained minimization formulation with AMEn for the global hybrid IAS with $r=-1$ and the 3D model problem for different levels of discretization.}  \label{table:3Dsweeps}
\end{minipage}
\end{table}

\section{Outlook}

We demonstrated that combining low rank tensor methods with sparsity inducing hyperprior models gives us a powerful method for the reconstruction of sparse initial states for large scale PDE systems. We tested different parameter setups for the local an global version of the hybrid IAS algorithm and compared the performance of different local solvers. The IAS algorithm requires a number of subsequent solutions of a large scale minimization problem. We can solve the arising problems with standard methods for quadratic minimization like CLGS or formulate a constrained optimization problem. In the context of  IGA discretization, both representations can be solved efficiently using tensor train calculations and the Alternating Minimal Energy solver (AMEn) designed to work with low rank Kronecker product based problems.

In combination with an Isogeometric space discretization, tensor train calculations can be very efficient. In future research, we want to enhance the methods performance by equipping it with suitable preconditioning and a more fitting local solver. This will allow us to solve sparse reconstruction problems for large-scale three dimensional discretizations. 

\section*{Acknowledgements}
The authors would like to thank Daniela Calvetti for her helpful insights and sharing of implementations.

The work of the authors was supported by the German Science Foundation (DFG) through grant 1742243256 - TRR 9.

\bibliographystyle{siam}
\bibliography{literature}
\end{document}